\documentclass[12pt]{amsart}
\usepackage{mathrsfs}
\usepackage{amssymb}
\usepackage{cite}
\usepackage[all]{xy}
\usepackage{graphicx}
\usepackage{comment}
\usepackage{textcomp}
\usepackage{charter}
\usepackage[vflt]{floatflt}
\usepackage{mathtools}
\usepackage{enumerate}
\usepackage{wasysym}
\usepackage{rotating}
\usepackage{pdflscape}
\usepackage{fullpage}
\usepackage{bm}
\usepackage[pdftex,colorlinks]{hyperref}
\hypersetup{pdfauthor={Shamovich, E.},
bookmarksnumbered,
pdfstartview={FitH}}

\newtheorem{thm}{Theorem}[section]
\newtheorem{prop}[thm]{Proposition}
\newtheorem{lem}[thm]{Lemma}
\newtheorem{cor}[thm]{Corollary}
\newtheorem*{thm*}{Theorem}
\newtheorem*{thmA}{Theorem A}
\newtheorem*{thmB}{Theorem B}
\newtheorem*{cor*}{Corollary}

\theoremstyle{definition}
\newtheorem{dfn}[thm]{Definition}

\theoremstyle{remark}
\newtheorem{remark}[thm]{Remark}

\numberwithin{equation}{section}

\newcommand{\PGL}{\mathrm{PGL}}
\newcommand{\SL}{\mathrm{SL}}
\newcommand{\PU}{\mathrm{PU}}
\newcommand{\B}{\mathbb{B}}
\DeclareMathOperator{\sat}{sat}
\DeclareMathOperator{\row}{row}
\DeclareMathOperator{\qinn}{q-Inn}
\DeclareMathOperator{\aut}{Aut}

\title{Similarity preserving noncommutative maps}

\author{Eli Shamovich}
\address{Department of Mathematics\\
 Ben-Gurion University of the Negev\\
 Beer-Seva\; 8410501\\
 Israel}
\thanks{E.S. was partially supported by BSF grant 2022235}

\usepackage{mathrsfs}
\newcommand{\Hom}{\operatorname{Hom}} %External Hom
\newcommand{\depth}{\operatorname{depth}}

\newcommand{\Ext}{\operatorname{Ext}} %Global Ext
\newcommand{\GL}{\operatorname{{\mathbf GL}}}

\newcommand{\tr}{\operatorname{tr}}
\newcommand{\cE}{{\mathcal E}}

\newcommand{\cH}{{\mathcal H}}

\newcommand{\cK}{{\mathcal K}}
\newcommand{\cL}{{\mathcal L}}

\newcommand{\cO}{{\mathcal O}}
\newcommand{\cP}{{\mathcal P}}

\newcommand{\fB}{{\mathfrak B}}

\newcommand{\fm}{{\mathfrak m}}

\newcommand{\C}{{\mathbb C}}

\newcommand{\N}{{\mathbb N}}
\newcommand{\Z}{{\mathbb Z}}
\newcommand{\M}{{\mathbb M}}

\newcommand{\D}{{\mathbb D}}

\newcommand{\bH}{{\mathbb H}}

\newcommand{\bB}{\ensuremath{{\mathbf B}}}

\begin{document}

\begin{abstract}
    In this note we study similarity preserving automorphisms of noncommutative (NC for short) balls. We show that every similarity preserving NC automorphism is, in fact, given by a similarity with an NC invertible function. As an application, we show that the group of automorphisms of the free semigroup algebra is generated by unitarily implemented and inner automorphisms.
\end{abstract}
\maketitle

\section{Introduction}

This paper grew out of the joint work of the author with Salomon and Shalit \cite{SSS1,SSS2}. In these works we have studied the isomorphism problem for algebras that are quotients of the free semigroup algebra $\bH^{\infty}_d$ by WOT closed ideals. The free semigroup algebra is the WOT closed algebra generated by the creation operators on the full Fock space (see Section \ref{sec:quasi-inner} for more details). This algebra was extensively studied by many (see, for example, \cite{DavidsonPitts-inv, DavidsonPitts-alg, DavidsonPitts-NPint, Popescu-aut, JuryMartin-affil ,Popescu-funcI, Popescu-funcII, Popescu-func_calc}). These algebras appear as noncommutative generalizations of the Banach algebra of bounded analytic functions on the unit disc both in dilation theory \cite{Bunce, Frazho, Popescu-dilations} and in noncommutative analysis \cite{AgMcC-freepick,BMV-ncrkhs, BMV-int, JMS-bso, JMS-rat, JMS-clark, Popescu-aut, Popescu-funcI}. The particular result of interest for this paper is a similarity invariant version of Cartan's uniqueness theorem \cite[Theorem 4.15]{SSS2} and its conclusion for a group of automorphisms of $\bH^{\infty}_d$ \cite[Corollary 6.8]{SSS2}. The latter group is defined by Davidson and Pitts \cite{DavidsonPitts-alg} as the group of automorphisms that induce the identity map modulo the commutator ideal. Recall that Cartan's classical uniqueness theorem says that if $U \subset \C^d$ is a bounded domain and $f \colon U \to U$ is analytic such that for some $z_0 \in U$, $f(z_0) = z_0$ and $Df(z_0) = I_d$, then $f$ is the identity. In the noncommutative setting, \cite[Theorem 4.15]{SSS2} says that if $f$ is a noncommutative analytic self-map of the similarity envelope of the unit row ball that fixes $0$ and whose derivative at $0$ is the identity, then $f$ must preserve similarity orbits of irreducible points. Intuitively speaking, $f$ is the identity modulo similarity. It is a natural question to ask whether such a map must be implemented by a similarity with an analytic noncommutative function. Our main result answers this question positively (we refer the reader to the following section for the definitions).
\begin{thmA}
Let $\cE$ be an operator space structure on $\C^d$ for $d \geq 2$. Let $\bB_{\cE}$ be the open unit ball of $\cE$. Let $\widetilde{\bB_{\cE}}$ be the similarity envelope of $\bB_{\cE}$. If $\varphi \colon \widetilde{\bB_{\cE}} \to \widetilde{\bB_{\cE}}$ is a noncommutative analytic biholomorphism that preserves the similarity orbits of irreducible points, then there exists a noncommutative analytic invertible function $G$ such that $\varphi(X) = G(X)^{-1} X G(X)$ for all $X \in \widetilde{\bB_{\cE}}$.
\end{thmA}
This theorem can be viewed as a strengthening of \cite[Theorem 4.15]{SSS2} since it provides a description of all such automorphisms. Essentially, they are the identity up to similarity with an invertible analytic function. We then use this result to prove the following theorem.
\begin{thmB}
    Every quasi-inner automorphism of $\bH^{\infty}_d$ is inner.
\end{thmB}
This is a somewhat surprising (to the author) result, since it implies that the automorphism group of $\bH^{\infty}_d$ is generated by unitarily implemented automorphisms that arise from biholomorphic automorphisms of the Euclidean open unit ball in $\C^d$ and inner automorphisms.

\subsection*{Structure of the paper} The paper has four sections apart from the current one. Section \ref{sec:background} is dedicated to collecting the necessary information about noncommutative functions and domains. Section \ref{sec:good_dom} contains preliminary results on when the levels of the similarity envelopes of certain noncommutative domains are domains of holomorphy (equivalently, Stein domains). Section \ref{sec:similarity} is the technical heart of the paper. It contains the proof of Theorem A. Lastly, the necessary information about the free semigroup algebra and its automorphisms is contained in Section \ref{sec:quasi-inner}. In particular, Theorem B is proved in Section \ref{sec:quasi-inner}.

\subsection*{Tool usage} This paper consists of old ideas of the author developed during the author's stay at the University of Waterloo. Lemmas \ref{lem:direct_sum_cocycle} and \ref{lem:respect_direct_sum} were the missing ingredient and were provided by ChatGPT 6.0 Astra. Additionally, ChatGPT 6.0 Astra provided the key observation in the proof of Theorem \ref{thm:qinn_is_inn}. The paper itself was written by the author and ChatGPT 6.0 Astra was used to verify some of the propositions and check the spelling.

\section{Background and notations} \label{sec:background}

In this section, we recall some basics from the theory of noncommutative analysis. The interested reader is referred to the books \cite{KVV-book, AMcCY-book} for details. Let us fix $d \in \N$ and consider the free algebra $\C \langle z_1,\ldots,z_d \rangle$. For a finite word $\alpha$ over the alphabet $\{1,\ldots,d\}$, we will write $|\alpha|$ for the length of $\alpha$ and $z^{\alpha} = z_{\alpha_1} \cdots z_{\alpha_n}$, where $|\alpha| = n$. We will call the elements of the free algebra noncommutative (NC for short) polynomials or simply polynomials when there can be no confusion. The finite-dimensional representations of the free algebra are parametrized by $d$-tuples of matrices. We denote this space by
\[
\M_d = \bigsqcup_{n=1}^{\infty} M_n(\C)^d.
\]
In the NC setting, $\M_d$ is analogous to $\C^d$, where the latter is the space of one-dimensional representations of $\C[z_1,\ldots,z_d]$. More generally, given an operator space $\cE$, we can consider $\M(\cE) = \bigsqcup_{n=1}^{\infty} M_n(\cE)$. If we fix an operator space structure on $\C^d$ and identify $M_n(\C^d) \cong M_n(\C)^d$, then $\M_d = \M(\cE)$. The topology on each level is independent of the choice of the operator structure due to the fact that all operator space structures on $\C^d$ are completely boundedly isomorphic.

The space $\M_d$ comes equipped with two natural operations arising from the direct sum and isomorphism of representations. For $X \in M_n(\C)^d$ and $Y \in M_m(\C)^d$, we define the \emph{direct sum} by
\[
X \oplus Y = \left( \begin{pmatrix} X_1 & 0 \\ 0 & Y_1 \end{pmatrix},\, \ldots,\, \begin{pmatrix} X_d & 0 \\ 0 & Y_d \end{pmatrix}\right) \in M_{n+m}(\C)^d.
\]
The direct sum is used to connect the different strata of $\M_d$. For \emph{similarity}, we set for $S \in \GL_n(\C)$ and $X \in M_n(\C)^d$,
\[
S X S^{-1} = \left( S X_1 S^{-1},\, \ldots, S X_d S^{-1} \right).
\]
Clearly, every $f \in \C \langle z_1,\ldots,z_d \rangle$ satisfies the following three properties:
\begin{itemize}
    \item[(i)] $f$ is \emph{graded}, i.e., $f(M_n(\C)^d) \subset M_n(\C)$;

    \item[(ii)] $f$ \emph{respects direct sums}, i.e., for every $X, Y \in \M_d$, $f(X \oplus Y) = f(X) \oplus f(Y)$;

    \item[(iii)] $f$ \emph{respects similarities}, i.e., $f(S X S^{-1}) = S f(X) S^{-1}$.
\end{itemize}
We say that a point $X \in \M_d$ is \emph{irreducible} if the representation of $\C \langle z_1,\ldots,z_d \rangle$ induced by $X$ is irreducible. In other words, the algebra generated by the $X_j$ is all of $M_n(\C)$. By Burnside's theorem, this is equivalent to the fact that the $X_j$ have no non-trivial common invariant subspace. Otherwise, we will say that $X$ is \emph{reducible}. If $X \in M_n(\C)^d$ is reducible, there exists a unitary $U$ such that
\[
U X U^* = \begin{pmatrix} X^{(1)} & \star & \star & \star \\ 0 & X^{(2)} & \star & \star \\ 0 & 0 & \ddots & \star \\ 0 & \cdots & 0 & X^{(k)} \end{pmatrix}.
\]
Here the points $X^{(j)}$ are irreducible. We will call the  similarity orbits of the $X^{(j)}$ the Jordan-H\"{o}lder constituents of $X$. We write $X^{ss} = \oplus_{j=1}^k X^{(j)}$. The latter point is defined up to a similarity. We call $X^{ss}$ the \emph{semisimplification} of $X$.

Since we want to deal with the noncommutative analog of analytic functions, we need to consider subdomains of $\M_d$. We say that $\Omega \subseteq \M_d$ is an NC set if $\Omega$ is closed under direct sums. We write $\Omega_n = \Omega \cap M_n(\C)^d$. We say that an NC set $\Omega$ is an NC domain if for every $n \in \N$, $\Omega_n$ is an open connected set that is invariant under similarities with unitary matrices. Lastly, we denote by $\widetilde{\Omega}$ the \emph{similarity envelope} of $\Omega$. Namely,
\[
\widetilde{\Omega} = \bigsqcup_{n=1}^{\infty} \left\{ S X S^{-1} \mid S \in \GL_n(\C) \text{ and } X \in \Omega_n\right\}.
\]
A class of NC domains central to this paper consists of the open balls of operator spaces. Let $\cE$ be an operator space structure on $\C^d$. We set
\[
\bB_{\cE} = \bigsqcup_{n=1}^{\infty} \left\{ X \in M_n(\C)^d \mid \|X\|_{\cE,n} < 1\right\}.
\]
In particular, the unit ball of the row operator space is
\[
\fB_d = \bigsqcup_{n=1}^{\infty} \left\{ X \in M_n(\C)^d \mid \sum_{j=1}^d X_j X_j^* < 1\right\}.
\]
The operator space balls are a particular example of matrix convex sets (see, for example, \cite{DavidsonKennedy_choquet,HKM12}). An NC set is matrix convex if it is closed under compression by isometries. Namely, if $X \in \Omega_n$ and $V \colon \C^m \to \C^n$ is an isometry, then
\[
V^*X V = \left( V^*X_1V,\, \ldots,\, V^*X_d V\right) \in \Omega_m.
\]
A more general class of NC domains was studied in \cite{AgMcC-implicit, AgMcC-freepick} in the case of polynomial maps and in \cite{BMV-int,BMV-ncrkhs} in general. Namely, let $\Omega \subseteq \M_d$ be an NC domain and let $Q \colon \Omega \to \M(B(\cH,\cK))$ be an analytic NC function, where $\cH$ and $\cK$ are separable Hilbert spaces. Then, we set
\[
D_Q = \left\{ X \in \Omega \mid \|Q(X)\| < 1 \right\}.
\]
Using the completely isometric embedding of an operator space into $B(\cH)$ we can see that every ball of an operator space is a domain of the type $\D_Q$. In \cite{BMV-ncrkhs}, Ball, Marx, and Vinnikov used the Effros-Winkler Hahn-Banach theorem \cite{EffrosWinkler} to show that every matrix convex set is of the form $D_Q$. In general, a set of this form need not have every level connected. Hence, we will make the implicit assumption that $\D_Q$ is connected, whenever we discuss such a domain. Our main focus in the paper will be on operator space balls. Therefore, this assumption will not play a big role. We will also make the assumption that $\Omega$ is the maximal domain for $Q$. Namely, there is no NC domain $\Omega \subsetneq \Omega' \subseteq \M_d$ and an NC analytic function $Q' \colon \Omega' \to \M(B(\cH,\cK))$ such that $Q'|_{\Omega} = Q$.

Given an open NC domain $\Omega \subset \M_d$, a function $f \colon \Omega \to \M_1$ is called an NC function if it satisfies the above properties (i)-(iii) relative to $\Omega$. A surprising fact is that if an NC function is locally bounded on $\Omega$, then it is automatically levelwise analytic \cite{AgMcC-global,KVV-book}. By the appendix in \cite{KVV-book}, every NC function $f$ extends uniquely to $\widetilde{\Omega}$. Moreover, every function admits an expansion around every point in $\Omega$ into a locally uniformly and absolutely convergent ``power series'' \cite{KVV-book}. If the point is scalar, then the expansion is a power series in free variables (i.e., an element of the formal power series ring $\C\langle \langle z_1,\ldots, z_d \rangle \rangle$).  

Let $U \subseteq M_n(\C)^d$ be a similarity invariant open set and let $f \colon U \to M_n(\C)$ be a holomorphic function. Following \cite{GMS-boundaries}, we say that $f$ is a \emph{concomitant} if $f$ is $\GL_n(\C)$-equivariant.

\section{Good NC domains} \label{sec:good_dom}

The goal of this section is to study a class of NC domains that are suitable as representation spaces of operator algebras. Recall that for an operator algebra $A \subset B(\cH)$, the restrictions to invariant subspaces and compressions to co-invariant subspaces are completely contractive representations. Moreover, every bounded representation of $A$ into $M_n$ is automatically completely bounded and by Paulsen's similarity theorem is similar to a completely contractive representation. This leads us to the following definition.
\begin{dfn}
    Let $\Omega \subset \M_d$ be an NC domain. We say that $\Omega$ is \emph{good} if it is closed under compressions to invariant and co-invariant subspaces.
\end{dfn}
To be precise, if $\Omega$ is a good domain and $X \in \Omega_n$, we say that $V \subset \C^n$ is invariant if $V$ is invariant under the coordinates of $X$. Hence, if $\iota \colon \C^k \to \C^n$ is an isometry with range $V$, then $\iota^* X \iota \in \Omega_k$. The same holds for coinvariant subspaces. If $\Omega$ is matrix convex, then it is trivially good. 

Since every NC domain $\Omega$ is invariant under unitary similarities, if $X$ has an invariant subspace, then we can find a unitary $U$ such that
\[
U X U^* = \begin{pmatrix} X_1 & \star \\ 0 & X_2 \end{pmatrix}.
\]
Hence, $X_1, X_2 \in \Omega$ if and only if $\Omega$ is good. This leads us to the following lemma.
\begin{lem} \label{lem:d_q_good}
    Let $\Omega \subseteq \M_d$ be a good NC domain and $Q \colon \Omega \to \M(B(\cH,\cK))$ an NC function. Then, $\D_Q$ is a good NC domain if every level is connected.
\end{lem}
\begin{proof}
    Since $Q$ respects unitary similarities, for every $X \in \D_{Q,n}$ with an invariant subspace, we have that
    \[
    U X U^* = \begin{pmatrix} X_1 & \star \\ 0 & X_2 \end{pmatrix},\, Q(X) = (I \otimes U^*) \begin{pmatrix} Q(X_1) & \star \\ 0 & Q(X_2) \end{pmatrix}(I \otimes U).
    \]
    Since $X_1, X_2 \in \Omega$, we have that $\|Q(X_1)\|, \|Q(X_2)\| \leq \|Q(X)\| < 1$. We conclude that $\D_Q$ is good.
\end{proof}

\begin{remark}
    The proof works even if the levels are not necessarily connected.
\end{remark}

The following lemma is immediate.
\begin{lem} \label{lem:good_jh}
    Let $\Omega \subseteq \M_d$ be a good domain. Then for every $X \in \Omega$, a representative of every Jordan-H\"{o}lder constituent of $X$ is contained in $\Omega$. Therefore, a representative of the semisimplification of $X$ is in $\Omega$.
\end{lem}

Let us denote by $M_n(\C)^d//\GL_n(\C)$ the categorical quotient by the similarity action. Namely, we define a relation on $M_n(\C)^d$ by $X \sim Y$ if $\overline{\GL_n(\C) \cdot X} \cap \overline{\GL_n(\C) \cdot Y} \neq \emptyset$. The set $M_n(\C)^d//\GL_n(\C)$ is the quotient by this equivalence relation and we write $p_n$ for the quotient map. We collect a few facts about the relation and $p_n$ that we will require. The references for these facts are \cite[Section 12.6]{Artin-azumaya} and \cite{Procesi-invariants}. Firstly, for $X \in \M_d$ irreducible, $p_n^{-1}(p_n(X))$ is the similarity orbit of $X$. More generally, for a reducible $X$, the fiber $p_n^{-1}(p_n(X))$ contains exactly one closed orbit, namely the orbit of the semisimplification of $X$. In particular, $X \sim Y$ if and only if $X^{ss} \sim Y^{ss}$. We say that $E \subset M_n(\C)^d$ is \emph{saturated} if $p_n^{-1}(p_n(E)) = E$.
\begin{dfn}
    Let $\Omega \subset \M_d$. We define the saturation of $\Omega$ by
    \[
    \sat(\Omega) = \bigsqcup_{n=1}^{\infty} p_n^{-1}(p_n(\Omega_n)).
    \]
    We say that $\Omega$ is \emph{saturated} if $\Omega = \sat(\Omega)$.
\end{dfn}
We record a lemma that is immediate from the above definition.
\begin{lem} \label{lem:saturatation_properties}
    Let $\Omega \subseteq \M_d$. Then,
    \begin{itemize}
        \item[(i)] $\sat(\sat(\Omega)) = \sat(\Omega)$. 

        \item[(ii)] $\widetilde{\Omega} \subset \sat(\Omega) = \sat(\widetilde{\Omega})$.
    \end{itemize}
\end{lem}
The next proposition ties the notion of a good NC domain to saturation.
\begin{prop} \label{prop:good_saturated}
    Let $\Omega \subseteq \M_d$ be a good NC domain. Then, $\widetilde{\Omega}$ is also a good domain. Moreover, $\widetilde{\Omega}$ is saturated.
\end{prop}
\begin{proof}
    For the first claim, we consider $X \in \Omega_n$. Let us fix $V \subset \C^n$ a non-trivial $X$-invariant subspace. Let $S \in \GL_n(\C)$. Then, $S V$ is invariant for $S X S^{-1}$. Let $\alpha, \beta \colon \C^k \to \C^n$ be isometries with ranges $V$ and $S V$, respectively. Then, $\alpha^* S^{-1} \beta \in \GL_k(\C)$. Moreover, let $\xi \in \C^k$. Then there exists $\eta \in \C^k$ such that $\beta \xi = S \alpha \eta$. Hence,
    \[
    (\beta^* S \alpha)(\alpha^* S^{-1} \beta) \xi = \beta^* S \alpha \eta = \xi.
    \]
    Now since $\Omega$ is good, $\alpha^* X \alpha \in \Omega_k$. Now we use the fact that $V$ is $X$-invariant and compute
    \[
    (\beta^* S \alpha)(\alpha^* X_j \alpha)(\alpha^* S^{-1}\beta) \xi = \beta^* S X_j \alpha \alpha^* S^{-1} \beta \xi = \beta^* S X_j \alpha \eta = \beta^* S X_j S^{-1} \beta \xi.
    \]
    We conclude that $\beta^* S X S^{-1} \beta \in \widetilde{\Omega}$. A similar argument works for coinvariant spaces.
    
    Now we prove the second claim. By the preceding lemma, it suffices to show that $\sat(\Omega) \subseteq \widetilde{\Omega}$. Now let $X \in \Omega_n$. If $X$ is irreducible, then $p_n^{-1}(p_n(X)) \subseteq \widetilde{\Omega}$. If $X$ is reducible, then $X^{ss}$ is contained in $\Omega$, since $\Omega$ is good. Since $p_n^{-1}(p_n(X))$ is determined by the semisimplifications, we have that for $Y \in p_n^{-1}(p_n(X))$, there exists $S \in \GL_n(\C)$ such that
    \[
    S Y S^{-1} = \begin{pmatrix} X^{(1)} & \star & \star & \star \\ 0 & X^{(2)} & \star & \star \\ 0 & 0 & \ddots & \star \\ 0 & \cdots & 0 & X^{(k)} \end{pmatrix}.
    \]
    For every $t >0$ define
    \[
    R_t = t^{k-1} I \oplus t^{k-2} I \oplus \cdots \oplus t I \oplus I.
    \]
    Then, it is not hard to check that $\lim_{t \to 0+} R_t S Y S^{-1} R_t^{-1} = X^{(1)} \oplus \cdots \oplus X^{(k)} \in \Omega_n$. Since $\Omega_n$ is open, we see that $Y \in \widetilde{\Omega}$ and we are done.
\end{proof}

\begin{cor} \label{cor:D_Q-good+sat}
    Let $\Omega \subseteq \M_d$ be a good NC domain and let $Q \colon \Omega \to \M(B(\cH,\cK))$ be an analytic NC function. Then, $\widetilde{\D_Q}$ is good and saturated.
\end{cor}

In \cite{SerKha-domains}, Sergeev and Kha{\u\i}ntsner studied the question of whether the orbit of the second level of the NC polydisc under the action of the group $\SL_2(\C) \times \SL_2(\C)$ is a domain of holomorphy. The answer to the above question is positive. To do this, they have coined the notion of orbital convexity. In \cite{SerKha-domains}, they study the action of a compact Lie group $K$ on the affine space $\C^d$ and the extension of the action to the complexification $K^{\C}$ of $K$. We restrict our attention to the case when $K = U_n$, the affine space is $M_n(\C)^d$, and the complexification is $K^{\C} = \GL_n(\C)$. Since our domains are NC, we consider the actions on all levels. The following definition is \cite[Definition in 2.3]{SerKha-domains} adapted to our setting.
\begin{dfn}
    Let $\Omega \subseteq \M_d$ be a good NC domain. We say that $\Omega$ is \emph{orbit convex} if for every $X \in \Omega_n$ and every selfadjoint $H \in M_n(\C)$ such that $e^{H} X e^{-H} \in \Omega_n$, we have that for every $0 \leq t \leq 1$, $e^{t H} X e^{- t H} \in \Omega_n$.
\end{dfn}
The preceding definition is equivalent to asking that for every $X \in \Omega$ and every strictly positive $A$ such that $A X A^{-1} \in \Omega_n$, we have that $A^t X A^{-t} \in \Omega_n$ for every $0 \leq t \leq 1$. We will use this observation in the proof of the following lemma.

\begin{lem} \label{lem:D_Q-hol_conv}
    Let $Q \colon \M_d \to \M(B(\cH,\cK))$ be an analytic NC function that does not extend to a larger NC domain. Then, for every $n \in \N$ such that $\D_{Q,n} \neq \emptyset$, $\widetilde{D_{Q,n}}$ is holomorphically convex. 
\end{lem}
\begin{proof}
    Our goal is to apply \cite[Theorem 3]{SerKha-domains}. To this end, we observe that by Corollary \ref{cor:D_Q-good+sat}, $\widetilde{\D_Q}$ is good and saturated. Let us fix $n \in \N$ such that $\D_{Q,n} \neq \emptyset$. Then, $\widetilde{\D_{Q,n}}$ is saturated. Moreover, $\D_{Q,n}$ is holomorphically convex. Indeed, if we write $\rho(X) = \|X\|_{\row} - \log(1 - \|Q(X)\|)$, then $\rho$ is a plurisubharmonic exhaustion function for $\D_{Q,n}$. Therefore, we need only to show that $\D_Q$ is orbit convex. To this end, let $X \in \D_{Q,n}$ and $A$ be strictly positive such that $A^{-1} X A \in \D_{Q,n}$. We need to show that $A^{-t} X A^{t} \in \D_{Q,n}$ for every $0 \leq t \leq 1$. Our assumption is equivalent to $Q(X) A^2 Q(X)^* < A^2$. We need to show that for every $0 \leq t \leq 1$, $Q(X) A^{2t} Q(X)^* < A^{2t}$. Now recall from \cite[Theorem V.1.9]{Bhatia-book} and \cite[Theorem V.2.5]{Bhatia-book} that for every $0 \leq t \leq 1$, the function $f_t(x) = x^t$ is operator monotone and concave on $[0,\infty)$. Therefore, applying \cite[Theorem V.2.3(iii)]{Bhatia-book} to $-f_t$, we have that $Q(X) A^{2t} Q(X)^* \leq A^{2t}$. Now if $T$ and $S$ are positive operators such that $T < S$, then there is an $\varepsilon > 0$ such that $T \leq S - \varepsilon$. Hence, by monotonicity and functional calculus, $T^t \leq (S -\varepsilon)^t < S^t$. Hence, we get the desired result.
\end{proof}

\begin{remark}
    One can show that this holds more generally. Namely, if $\Omega$ is a similarity invariant NC domain such that each non-empty level is a domain of holomorphy and $Q \colon \Omega \to \M(B(\cH,\cK))$ is an NC analytic function that does not extend analytically to any NC domain that properly contains $\Omega$, then the above lemma still holds. Since this will not be relevant to the main results of the paper, we do not proceed to work this out.
\end{remark}

\section{Similarity preservers} \label{sec:similarity}

In this section we will restrict ourselves to operator space balls. Let us fix an operator space structure $\cE$ on $\C^d$ and let $\bB_{\cE} \subseteq \M_d$ be the open unit ball of $\cE$. As mentioned in Section \ref{sec:background}, $\bB_{\cE} = \D_Q$ for a linear $Q$. Therefore, the results of the preceding section apply to $\bB_{\cE}$. In particular, we have that $\bB_{\cE}$ is good and orbit convex and $\widetilde{\bB_{\cE}}$ is saturated and each level is a domain of holomorphy. 

For the purposes of this section, we fix $\varphi \colon \widetilde{\bB_\cE} \to \widetilde{\bB_{\cE}}$ an NC analytic biholomorphism such that $\varphi(X) \in \widetilde{\{X\}}$ for every $X \in \bB_{\cE}$. Our goal is to show that $\varphi$ is implemented by a similarity with an invertible analytic NC function. Our strategy is to first prove that it is true on every sufficiently high level, and then extend it to lower levels and fix the direct sum property. Therefore, fix $n \geq 3$ and assume that $d \geq 2$. These assumptions ensure that the collection of all reducible points forms an algebraic subvariety of $M_n(\C)^d$ of codimension strictly greater than $1$.

Let us write $G_{d,n} \subset M_n(\C)^d$ for the set of all irreducible tuples and $E_{d,n} = (\widetilde{\bB_{\cE}})_n \cap G_{d,n}$. Furthermore, recall that $p_n \colon M_n(\C)^d \to M_n(\C)^d//\GL_n(\C)$ is the quotient map by the relation described in the preceding section. We denote by $V_{d,n} = p_n(G_{d,n})$ and $Q_{d,n} = p_n(E_{d,n})$. By the result of Le Bruyn and Procesi \cite[Theorem II.3.4]{LeBruynProcesi-etale}, $p_n \colon G_{d,n} \to V_{d,n}$ is a principal $\PGL_n(\C)$ bundle and, moreover, $V_{d,n}$ is the set of smooth points in $M_n(\C)^d//\GL_n(\C)$. 
\begin{lem} \label{lem:phi_induced}
    For every $X \in E_{d,n}$, there exists a unique $S(X) \in \PGL_n(\C)$ such that $\varphi(X) = S(X)^{-1} X S(X)$. Moreover, the function $S \colon E_{d,n} \to \PGL_n(\C)$ is holomorphic and $\PGL_n(\C)$-equivariant, i.e., $S(T X T^{-1}) = T S(X) T^{-1}$.
\end{lem}
\begin{proof}
    By assumption such an $S(X)$ exists for every $X \in \bB_{\cE,n}$. However, since $X$ is irreducible, the result of Procesi and the fact that $\varphi$ is bijective imply that this $S(X)$ is unique. Since $\varphi$ is equivariant, we have that for every $X \in E_{d,n}$ and $T \in \PGL_n(\C)$,
    \[
    S^{-1}(T X T^{-1}) T X T^{-1} S(T X T^{-1}) = \varphi(T X T^{-1}) = T \varphi(X) T^{-1} = T S^{-1}(X) X S(X) T^{-1}.
    \]
    By uniqueness, $S(T X T^{-1}) = T S(X) T^{-1}$. 

    Lastly, for $X \in E_{d,n}$ we consider the equation $X S = S \varphi(X)$. This is a linear equation with holomorphic coefficients. Moreover, the rank is constant on the open set $E_{d,n}$, since the commutant of $X$ is just the scalars. Therefore, there is a local holomorphic section, which we project onto $\PGL_n(\C)$ to obtain $S(X)$. Therefore, $S(X)$ is holomorphic.
    \end{proof}

    Following \cite{GMS-boundaries}, we set $L_{d,n} = E_{d,n} \times^{\PGL_n(\C)} \GL_n(\C)$ and $P_{d,n} = E_{d,n} \times^{\PGL_n(\C)} \PGL_n(\C)$. This is a standard notation for the quotient of the Cartesian product by the equivalence relation $(X, T) \sim (S^{-1} X S, S^{-1} T S)$. Using the projection on the first coordinate, these are $\PGL_n(\C)$-bundles on $Q_{d,n}$. The above lemma implies that $S(X)$ induces a holomorphic section of $P_{d,n}$. Moreover, we have the exact sequence of $\PGL_n(\C)$-bundles
    \[
    \xymatrix{1 \ar[r] & E_{d,n} \times ^{\PGL_n(\C)} \C^{\times} \ar[r] & L_{d,n} \ar[r] & P_{d,n} \ar[r] & 1.}
    \]
    Here the action of $\PGL_n(\C)$ on $\C^{\times}$ is trivial and the second map is the embedding of the center. Clearly, the sheaf of holomorphic sections of $E_{d,n} \times ^{\PGL_n(\C)} \C^{\times}$ is the sheaf of nonvanishing holomorphic functions on $Q_{d,n}$. As is rather common, we will write $\cO_{Q_{d,n}}$ for the sheaf of holomorphic functions on $Q_{d,n}$ and $\cO^{\times}_{Q_{d,n}}$ for the sheaf of invertible holomorphic functions on $Q_{d,n}$. Hence, if we write $\cL_{d,n}$ and $\cP_{d,n}$ for the sheaves of holomorphic sections of $L_{d,n}$ and $P_{d,n}$, respectively, we get an exact sequence
    \begin{equation} \label{eq:grouP_sheaf_seq}
    \xymatrix{0 \ar[r] & \cO^{\times}_{Q_{d,n}} \ar[r] & \cL_{d,n} \ar[r] & \cP_{d,n} \ar[r] & 0}.
    \end{equation}
    Our goal now is to show that $S$, which we have identified with an element of $H^0(Q_{d,n}, \cP_{d,n})$, can be lifted to an element of $H^0(Q_{d,n}, \cL_{d,n})$. The latter by \cite{GMS-boundaries} will give us the desired holomorphic concomitant implementing $\varphi$. This leads us to study the cohomology of these sheaves more closely. Recall that the exponential exact sequence is
    \begin{equation} \label{eq:exp_exact}
    \xymatrix{0 \ar[r] & \Z \ar[r] & \cO_{Q_{d,n}} \ar[r]^{\exp} & \cO^{\times}_{Q_{d,n}} \ar[r] & 0.}
    \end{equation}
    Passing to the long exact sequence in cohomology from the exponential exact sequence and using the fact that the sheaf cohomology of the constant sheaf on a manifold $Q_{d,n}$ is in fact the topological cohomology with integral coefficients, we see that we need to understand the topology of $Q_{d,n}$. This observation is useful, since the exact sequence \eqref{eq:grouP_sheaf_seq} is an exact sequence of group sheaves. Hence, by \cite[Section I.6, Theorem I.2]{Frenkel-nonab_coh}, there is an exact sequence in nonabelian cohomology
    \begin{equation} \label{eq:exact_seq_nonab}
    \xymatrix{0 \ar[r] & H^0(Q_{d,n},\cO_{Q_{d,n}}^{\times}) \ar[r] & H^0(Q_{d,n},\cL_{d,n}) \ar[r] & H^0(Q_{d,n},\cP_{d,n}) \ar[r] & H^1(Q_{d,n},\cO_{Q_{d,n}}^{\times}).}
    \end{equation}
    Therefore, to prove that we can lift $S$, it suffices to show that $H^1(Q_{d,n},\cO^{\times}_{Q_{d,n}}) = 0$. 
    \begin{lem} \label{lem:homotopy_E}
    Under our assumption that $d \geq 2$ and $n \geq 3$, we have that $E_{d,n}$ is path connected and $\pi_1(E_{d,n}) = \pi_2(E_{d,n}) = 0$.
    \end{lem}
    \begin{proof}
        Under our assumption, the complex codimension of the locus of reducible points is at least $2$. Therefore, its real codimension is at least $4$. Since $\widetilde{\bB_{\cE,n}}$ is contractible, $E_{d,n}$ is path connected. Now let $\gamma \colon S^1 \to E_{d,n}$ be a closed loop. Since $\widetilde{\bB_{\cE,n}}$ is contractible, there is an extension of $\gamma$ to the closed unit disc. By \cite[Theorem 2.2.2]{Hirsch-diff_top}, we may perturb our extension slightly and assume that it is smooth. Breaking up the singular locus into a union of smooth manifolds by iteratively taking the smooth points of this variety and then the smooth points of the remaining singular locus, we may apply the transversality theorem \cite[Theorem 3.2.5]{Hirsch-diff_top}. Indeed, the codimension is at least $4$ and the tangent space to the disc is two-dimensional. This argument shows that we can perturb our extension slightly to make it transverse to each stratum. Since our perturbations are arbitrarily small and the image of $\gamma$ is compact, we see that we can choose the perturbed extension to lie in $E_{d,n}$. We can now build a homotopy by moving $\gamma$ to the perturbed image of the circle and then contracting it. This proves that $\pi_1(E_{d,n}) = 0$. A similar argument works for $\pi_2(E_{d,n})$.
    \end{proof}

    The following lemma is well known, but we include it for the sake of completeness.
    \begin{lem} \label{lem:pgln_homotopy}
        For every $n \geq 2$, $\pi_1(\PGL_n(\C)) = \Z/n\Z$. 
    \end{lem}
    \begin{proof}
        By \cite[Theorem 14.3.11]{HilNeeb-lie_groups}, $\PGL_n(\C)$ deformation retracts onto its maximal compact subgroup $\PU_n$. Now we use the exact sequence,
        \[
        \xymatrix{1 \ar[r] & S^1 \ar[r] & \mathrm{U}_n \ar[r] & \PU_n \ar[r] & 1.}
        \]
        Using the long exact sequence in homotopy groups \cite[Theorem 4.41]{Hatcher-book} and the calculation in \cite[Example 4.55]{Hatcher-book}, we obtain
        \[
        \xymatrix{0 \ar[r] & \pi_2(\PU_n) \ar[r] & \Z \ar[r] & \Z \ar[r] & \pi_1(\PU_n) \ar[r] & 0.}
        \]
        Now the determinant induces an isomorphism $\pi_1(\mathrm{U}_n) \to \pi_1(S^1)$. Therefore, the fundamental group of the center sits inside $\pi_1(\mathrm{U}_n)$ as the subgroup of multiples of $n$. This observation combined with the exact sequence above proves the desired result.
    \end{proof}

    \begin{lem} \label{lem:homotopy_Q}
        The first homotopy groups of $Q_{d,n}$ are $\pi_1(Q_{d,n}) = 0$ and $\pi_2(Q_{d,n}) = \Z/n\Z$.
    \end{lem}
    \begin{proof}
        We apply the long exact sequence in homotopy groups of a fiber bundle to obtain
        \[
        \xymatrix{ \cdots \ar[r] & \pi_2(E_{d,n}) \ar[r] & \pi_2(Q_{d,n}) \ar[r] & \pi_1(\PGL_n(\C)) \ar[r] & \pi_1(E_{d,n}) \ar[r] & \pi_1(Q_{d,n}) \ar[r] & 0.}
        \]
        By Lemma \ref{lem:homotopy_E}, we get that $\pi_1(Q_{d,n}) = 0$ and $\pi_2(Q_{d,n}) \cong \pi_1(\PGL_n(\C)) = \Z/n\Z$, where the second equality follows from the preceding lemma.
    \end{proof}

    \begin{lem} \label{lem:cohomology_Q}
        We have that $H^1(Q_{d,n},\Z) = H^2(Q_{d,n},\Z) = 0$. Moreover, $H^1(Q_{d,n}, \cO_{Q_{d,n}}) = H^1(Q_{d,n},\cO^{\times}_{Q_{d,n}})$.
    \end{lem}
    \begin{proof}
        First, we note that by the Hurewicz theorem \cite[Theorem 4.32]{Hatcher-book}, $H_2(Q_{d,n},\Z) = \Z/n\Z$ and $H_1(Q_{d,n},\Z) = 0$. Now we apply the universal coefficient theorem \cite[Theorem 3.2]{Hatcher-book} to get the exact sequence for $i = 1,2$
        \[
        \xymatrix{ 0 \ar[r] & \Ext(H_{i-1}(Q_{d,n},\Z), \Z) \ar[r] & H^i(Q_{d,n},\Z) \ar[r] & \Hom(H_i(Q_{d,n}, \Z), \Z) \ar[r] & 0.}
        \]
        Now for $i = 2$, we immediately see that both the leftmost and rightmost terms are $0$. For $i=1$, the rightmost term is $0$ and $H_0(Q_{d,n},\Z) = \Z$ and, thus, the leftmost term is $0$.

        To see the last claim we use the long exact sequence in sheaf cohomology obtained from the exponential exact sequence \eqref{eq:exp_exact}. In particular, we get the exact piece of it
        \[
        \xymatrix{H^1(Q_{d,n},\Z) \ar[r] & H^1(Q_{d,n},\cO_{Q_{d,n}}) \ar[r] & H^1(Q_{d,n},\cO_{Q_{d,n}}^{\times}) \ar[r] & H^2(Q_{d,n}, \Z)}.
        \]
        Now we just apply the result on integral cohomology that we have obtained.
    \end{proof}

    The above lemma indicates that to show that our obstruction vanishes, it suffices to show that $H^1(Q_{d,n},\cO_{Q_{d,n}}) = 0$. To this end, we will use some local cohomology. We refer the reader to \cite[Chapter II]{BanSta-alg_methods} for details on local cohomology in the analytic category. Recall that a commutative local ring $(A,\fm)$ is called Cohen-Macaulay if $\depth_{\fm} A = \dim A$. Here the depth is the length of the longest regular sequence in $\fm$. We say that an analytic space $X$ is Cohen-Macaulay if all local rings are Cohen-Macaulay. We require the following lemma that summarizes some results on the relation between global and local sheaf cohomology.
    \begin{lem} \label{lem:locak_cohomology}
        Let $X$ be a pure-dimensional Cohen-Macaulay analytic space. Let $Y \subsetneq X$ be a closed analytic subspace of codimension at least $r \geq 1$. Then, for $0 \leq k < r - 1$, we have
        \[
        H^k(X,\cO_X) = H^k(X \setminus Y, \cO_{X \setminus Y}).
        \]
    \end{lem}
    \begin{proof}
        Assume that $\dim X = m$. Thus, $\dim Y \leq m - r$. Since $X$ is Cohen-Macaulay, we have that $\dim Y + (r-1) + 1 \leq m = \depth X$. By \cite[Corollary II.3.9]{BanSta-alg_methods}, for $q = r-1$, we get that the sheaves of local cohomology $\mathfrak{H}^k_Y = 0$ for $0 \leq k \leq r-1$. The equivalence of (c) and (d) in \cite[Theorem 3.6]{BanSta-alg_methods} implies the result.
    \end{proof}
    Since for every $n \in \N$, $\widetilde{\bB_{\cE,n}}$ is a domain of holomorphy, the above lemma combined with Cartan's theorem B \cite[Theorem 11.4.3]{Taylor-sev_comp_var} implies that for $d \geq 2$ and $n \geq 3$,  $H^0(E_{d,n},\cO_{E_{d,n}}) = H^0(\widetilde{\bB_{\cE,n}}, \cO_{\widetilde{\bB_{\cE,n}}})$ and for $1 \leq k < (d-1)(n-1) - 1$, $H^k(E_{d,n},\cO_{E_{d,n}}) = 0$. In particular, we can extend every holomorphic function on $E_{d,n}$ to $\widetilde{\bB_{\cE,n}}$ (this is in fact Hartogs' extension theorem in this case). To understand $Q_{d,n}$ we require another technical lemma.
    \begin{lem} \label{lem:CM}
        Let $X$ be a complex affine algebraic variety. If $X$ is Cohen-Macaulay as an algebraic variety, then it is Cohen-Macaulay as an analytic variety. Moreover, the dimensions of the algebraic and analytic local rings are the same.
    \end{lem}
    \begin{proof}
        This is \cite[Proposition XII.2.1]{SGA1}.
    \end{proof}

    \begin{lem} \label{lem:quotient_CM}
        The algebraic variety $M_n(\C)^d//\GL_n(\C)$ is an irreducible Cohen-Macaulay affine variety. Moreover, under the assumption that $d \geq 2$ and $n \geq 3$ the codimension of the singular locus inside $M_n(\C)^d//\GL_n(\C)$ is at least $3$.
    \end{lem}
    \begin{proof}
        Let us write $\C[M_n(\C)^d]$ for the (commutative) polynomial ring in the coordinates of the $d$-tuples of matrices. The first claim follows from the fact that $M_n(\C)^d//\GL_n(\C)$ is the spectrum of the ring $\C[M_n(\C)^d]^{\GL_n(\C)}$, which is Cohen-Macaulay by the main theorem of \cite{HochRob-CM}. The variety is irreducible as an image of an irreducible variety. The dimension of $M_n(\C)^d//\GL_n(\C)$ is $(d-1)n^2 + 1$, since the generic fiber of the quotient map is of dimension $n^2 - 1$. Now we note that the singular locus is the image of the reducible semisimple points. For every reducible semisimple point, we can break it into a direct sum. Hence, for $1 \leq m \leq n-1$, we have the maps
        \[
        M_m(\C)^d//\GL_m(\C) \times M_{n-m}(\C)^d//\GL_{n-m}(\C) \to M_n(\C)^d//\GL_n(\C).
        \]
        The images of these maps cover the singular locus. Hence, the dimension of the singular locus is bounded from above by $\max_{1 \leq m \leq n-1}\left((d-1)(m^2 + (n-m)^2) + 2\right) = (d-1)(1+(n-1)^2) + 2$. Therefore, the codimension of the singular locus is at least
        \[
        (d-1)n^2 + 1 - (d-1)(1 + (n-1)^2) - 2 = (d-1)(2n - 2) - 1 \geq 3.
        \]
         \end{proof}

    \begin{prop} \label{prop:levelwise_S}
        For every $d \geq 2$ and every $n \geq 3$, there is a matrix concomitant $S_n \colon \widetilde{\bB_{\cE,n}} \to \GL_n(\C)$ such that for every $X \in \widetilde{\bB_{\cE,n}}$, $\varphi(X) = S_n(X)^{-1} X S_n(X)$.
    \end{prop}
    \begin{proof}
        By \cite[Theorem 1]{SerKha-domains} and \cite[Corollary 2]{SerKha-domains}, we know that $p_n(\bB_{\cE,n})$ is a Stein space. Hence, by Cartan's theorem B, the first cohomology of its structure sheaf is $0$. By Lemmas \ref{lem:CM} and \ref{lem:quotient_CM}, we get that $H^1(Q_{d,n},\cO_{Q_{d,n}}) = 0$. Now Lemma \ref{lem:cohomology_Q} implies that $H^1(Q_{d,n},\cO_{Q_{d,n}}^{\times}) = 0$. Therefore, the exact sequence \eqref{eq:exact_seq_nonab} collapses to the short exact sequence
        \[
        \xymatrix{0 \ar[r] & H^0(Q_{d,n},\cO_{Q_{d,n}}^{\times}) \ar[r] & H^0(Q_{d,n},\cL_{d,n}) \ar[r] & H^0(Q_{d,n},\cP_{d,n}) \ar[r] & 0.}
        \]
        We conclude that we can lift $S$ to a concomitant on $E_{d,n}$. Now Lemma \ref{lem:locak_cohomology} and the paragraph following it imply that we can extend $S$ to a holomorphic matrix-valued function on $\widetilde{\bB_{\cE,n}}$. However, since we can extend $S^{-1}$ as well, the identity theorem implies that we obtain a concomitant $S \colon \widetilde{\bB_{\cE,n}} \to \GL_n(\C)$ such that for every $X \in \widetilde{\bB_{\cE,n}}$, $\varphi(X) = S(X)^{-1} X S(X)$.
    \end{proof}
    The exact sequence obtained in the proposition shows that $S$ is unique up to multiplication by a holomorphic $\GL_n(\C)$-invariant function. Let us fix on each level an invertible holomorphic concomitant $S_n$ that implements $\varphi$. To get the direct sum property, we need to find, for every $n \geq 3$, a holomorphic $\GL_n(\C)$-invariant function $\alpha_n$ such that $\alpha_n S_n$ respects direct sums. To this end we note that if $X \in \bB_{\cE,n}$ and $Y \in \bB_{\cE,m}$ for $n,m \geq 3$, then by the invariance property, we have that for every $\lambda, \mu \in \C^{\times}$
    \[
    (\lambda I_n \oplus \mu I_m) S_{n+m}(X \oplus Y) = S_{n+m}(X \oplus Y) (\lambda I_n \oplus \mu I_m).
    \]
    In particular, $S_{n+m}(X \oplus Y)$ is block-diagonal. Moreover,
    \begin{multline*}
    S_{n+m}(X \oplus Y)^{-1}(X \oplus Y) S_{n+m}(X \oplus Y) = \varphi(X \oplus Y) = \varphi(X) \oplus \varphi(Y) = \\  (S_n(X)^{-1} \oplus S_m(Y)^{-1}) (X \oplus Y) (S_n(X) \oplus S_m(Y)).
    \end{multline*}
    Write $S_{n+m}(X \oplus Y) = A(X, Y) \oplus B(X,Y)$. Then, $A(X,Y) S_n(X)^{-1}$ commutes with $X$. So for an irreducible $X$ it must be a non-zero scalar. This scalar is $a_{n,m}(X,Y) = \frac{1}{n}\tr\left(A(X,Y) S_n(X)^{-1}\right)$. Now fixing $Y$, we get a $\GL_n(\C)$-invariant function that we can extend as before to all of $\widetilde{\bB_{\cE,n}}$. The resulting function is nonvanishing on $E_{d,n}$. Any zeros would therefore lie in the reducible locus, whose complex codimension is at least $2$, whereas a nonempty zero locus of a nonzero holomorphic function has complex codimension $1$. The following lemma summarizes the properties of the above function. In particular, we show that our functions have a cocycle-like property.
    \begin{lem} \label{lem:direct_sum_cocycle}
    For every $n,m,k \geq 3$, we have the following two statements:
    \begin{itemize}
        \item[(i)] There exists a $\GL_n(\C) \times \GL_m(\C)$-equivariant function $a_{n,m} \colon \widetilde{\bB_{\cE,n}} \times \widetilde{\bB_{\cE,m}} \to \C^{\times}$ such that for every $(X,Y)$, we have that
        \[
        S_{n+m}(X \oplus Y) = \begin{pmatrix} a_{n,m}(X,Y) S_n(X) & 0 \\ 0 & a_{m,n}(Y,X) S_m(Y)\end{pmatrix}.
        \]
        \item[(ii)] For every $X \in \widetilde{\bB_{\cE,n}}$, $Y \in \widetilde{\bB_{\cE,m}}$, and $Z \in \widetilde{\bB_{\cE,k}}$, we have
        \[
        a_{n+m,k}(X \oplus Y, Z) a_{n,m}(X,Y) = a_{n,m+k}(X, Y \oplus Z).
        \]
    \end{itemize} 
    \end{lem}
    \begin{proof}
        The first part of the lemma, except for the identity, was proved in the discussion preceding it. To see the identity, we note that 
        \[
        B(X,Y) \oplus A(X,Y) = \begin{pmatrix} 0 & I_m \\ I_n & 0 \end{pmatrix} S_{n+m}(X \oplus Y) \begin{pmatrix} 0 & I_n \\ I_m & 0 \end{pmatrix} = S_{n+m}(Y \oplus X).
        \]
        Now using the same considerations for $S_{n+m}(Y \oplus X)$, we get the desired identity.

        For the second part, we use the identity from the first part of the lemma
        \begin{multline*}
        S_{n+m+k}(X \oplus Y \oplus Z) = \begin{pmatrix} a_{n+m,k}(X \oplus Y, Z) S_{n+m}(X \oplus Y) & 0 \\ 0 &  a_{k,n+m}(Z, X \oplus Y) S_k(Z) \end{pmatrix} = \\ \begin{pmatrix} a_{n+m,k}(X \oplus Y, Z) a_{n,m}(X,Y) S_n(X) & 0 & 0 \\ 0 & a_{n+m,k}(X \oplus Y, Z) a_{m,n}(Y,X) S_m(Y) & 0 \\ 0 & 0 & a_{k,n+m}(Z, X \oplus Y) S_k(Z) \end{pmatrix}.
        \end{multline*}
        On the other hand, we have that
        \[
        S_{n+m+k}(X \oplus Y \oplus Z) = \begin{pmatrix} a_{n,m+k}(X, Y \oplus Z) S_n(X) & 0 \\ 0 & a_{m+k,n}(Y \oplus Z, X) S_{m+k}(Y \oplus Z) \end{pmatrix}.
        \]
        Now comparing the upper left corners and using the invertibility of $S_n(X)$, we get the cocycle-like property.
    \end{proof}

    \begin{lem} \label{lem:respect_direct_sum}
        Fix $Z_0 \in \widetilde{\bB_{\cE,3}}$ and define for every $n \geq 3$, 
        \[
        \alpha_n(X) = \dfrac{a_{n,3}(X,Z_0)}{a_{3,n}(Z_0,X)} \text{ and } G_n(X) = \alpha_n(X) S_n(X).
        \]
        Then, $\alpha_n$ is a non-vanishing analytic $\GL_n(\C)$-equivariant function on $\widetilde{\bB_{\cE,n}}$ and for every $n, m \geq 3$ and $X \in \widetilde{\bB_{\cE,n}}$ and $Y \in \widetilde{\bB_{\cE,m}}$, $G_{n+m}(X \oplus Y) = G_n(X) \oplus G_m(Y)$.
    \end{lem}
    \begin{proof}
        First we observe that by the invariance property of $a_{n,m}$, we have that for all $n,m,k \geq 3$, $a_{n+m,k}(X \oplus Y,Z) = a_{n+m,k}(Y \oplus X,Z)$ and similarly for the second coordinate. The first claim about $\alpha_n$ is clear from the preceding lemma. Now let us fix $n,m \geq 3$ and compute
        \begin{multline*}
            G_{n+m}(X \oplus Y) = \alpha_{n+m}(X\oplus Y) S_{n+m}(X \oplus Y) = \\ \alpha_{n+m}(X\oplus Y) \begin{pmatrix} a_{n,m}(X,Y) S_n(X) & 0 \\ 0 & a_{m,n}(Y,X) S_m(Y)\end{pmatrix}.
        \end{multline*}
        Now
        \begin{multline*}
        \alpha_{n+m}(X\oplus Y) a_{n,m}(X,Y) = \dfrac{a_{n+m,3}(X \oplus Y,Z_0) a_{n,m}(X,Y)}{a_{3,n+m}(Z_0, X \oplus Y)} = \dfrac{a_{n,m+3}(X, Y \oplus Z_0)}{a_{3,n+m}(Z_0, X \oplus Y)} = \\ \dfrac{a_{n,m+3}(X, Z_0 \oplus Y)}{a_{3,n+m}(Z_0, X \oplus Y)} = \dfrac{a_{n+3,m}(X \oplus Z_0, Y) a_{n,3}(X,Z_0)}{a_{n+3,m}(Z_0 \oplus X,Y)a_{3,n}(Z_0,X)} = \alpha_n(X).
        \end{multline*}
        A similar argument works for the second summand.
    \end{proof}
    Now we are ready to prove the main result of the section.
    \begin{thm} \label{thm:main}
        Let $\cE$ be an operator space structure on $\C^d$ with $d \geq 2$ and let $\varphi \colon \widetilde{\bB_{\cE}} \to \widetilde{\bB_{\cE}}$ be an NC analytic biholomorphism such that for every $X \in \bB_{\cE}$, $\varphi(X) \in \widetilde{\{X\}}$. Then, there exists an NC analytic function $G$ on $\widetilde{\bB_{\cE}}$ such that for every $X \in \widetilde{\bB_{\cE}}$,
        \[
        \varphi(X) = G(X)^{-1} X G(X).
        \]
    \end{thm}
    \begin{proof}
        The preceding lemma allows us to construct an analytic NC function $G$ on levels $n \geq 3$ by setting $G(X) = G_n(X)$ for $X \in \widetilde{\bB_{\cE,n}}$. We need only to extend this function to the lower levels. To this end let $X \in \widetilde{\bB_{\cE,m}}$ for $m = 1$ or $m = 2$. Now consider $I_3 \otimes X = X^{\oplus 3}$. Then, by equivariance of $G$, for every permutation matrix $U$ of dimension $3$, we have that
        \[
        G(I_3 \otimes X) = G\bigl((U^* \otimes I_m)(I_3 \otimes X)(U \otimes I_m)\bigr) = (U^* \otimes I_m) G(I_3 \otimes X) (U \otimes I_m).
        \]
        Combining this with the same argument with $U$ replaced by a diagonal matrix, we get that $G(I_3 \otimes X) = I_3 \otimes G_m(X)$. Here $G_m(X)$ is automatically analytic and $\GL_m(\C)$-equivariant. Therefore, we have extended $G$ to the lower levels. Moreover, since $\varphi$ respects direct sums, we get that the identity $\varphi(X) = G(X)^{-1} X G(X)$ is satisfied on all levels now. If we take $k > 3$ and consider $I_k \otimes X$, then, as before, $G(I_k \otimes X) = I_k \otimes H(X)$. However, because of the fact that $G$ respects direct sums for $n \geq 3$, we have that
        \[
        I_{3k} \otimes G_m(X) = G(I_k \otimes I_3 \otimes X) = G(I_{3k} \otimes X) = G(I_3 \otimes I_k \otimes X) = I_{3k} \otimes H(X).
        \]
        Hence, the definition is independent of the choice of the size of the amplification. Now let $X, Y \in \widetilde{\bB_{\cE}}$. Then, $I_3 \otimes (X \oplus Y)$ is necessarily on a high enough level. Now taking a  shuffle unitary $P$ such that $P(I_3 \otimes (X \oplus Y))P^* = (I_3 \otimes X) \oplus (I_3 \otimes Y)$, we get that 
        \[
        I_3 \otimes G(X \oplus Y) = G(I_3 \otimes (X \oplus Y)) = P^* (G(I_3 \otimes X) \oplus G(I_3 \otimes Y)) P = I_3 \otimes (G(X) \oplus G(Y)).
        \]
        Hence, $G$ respects direct sums and we are done.
    \end{proof}

    \begin{cor} \label{cor:only_irreducible}
        Let $\varphi \colon \widetilde{\bB_{\cE}} \to \widetilde{\bB_{\cE}}$ be an NC analytic biholomorphism such that for every irreducible $X \in \bB_{\cE}$, $\varphi(X) \in \widetilde{\{X\}}$. Then, $\varphi$ preserves similarity orbits and is of the form $\varphi(X) = G(X)^{-1} X G(X)$ for some invertible NC analytic function $G$ on $\widetilde{\bB_{\cE}}$.
    \end{cor}
    \begin{proof}
        Since $\varphi$ preserves the similarity orbits of irreducible points, we have that for every $n \in \N$, $\varphi$ induces the identity map on $p_n(\bB_{\cE,n}) \subseteq M_n(\C)^d//\GL_n(\C)$. Moreover, the entire proof strategy consists of lifting the similarity on $E_{d,n}$ to an invertible equivariant function. This part translates verbatim. Now the identity theorem for analytic functions ensures that the equality $\varphi(X) = S_n(X)^{-1} X S_n(X)$ holds for $d \geq 2$ and $n \geq 3$ on the reducible locus, as well. The rest of the proof proceeds identically.
    \end{proof}

    \begin{remark} \label{rem:proof_extension}
        Note that the restriction to operator space balls is made to simplify the assumptions. Clearly, the proof technique extends to domains of the form $\D_Q$ that have sufficiently nice topology at all sufficiently high levels.
    \end{remark}

\section{Automorphisms of $\bH^{\infty}_d$} \label{sec:quasi-inner}

The full Fock space $\bH^2_d$ is the completion of the free algebra $\C \langle z_1,\ldots, z_d \rangle$ with respect to the norm induced by the inner product that makes the monomials an orthonormal basis. Left multiplication by the coordinates gives us a row isometry. We will write $L_{z_i} = L_i$ and $L = (L_1,\ldots,L_d) \colon \bH^2_d \otimes \C^d \to \bH^2_d$. The WOT closed algebra generated by the $L_i$ will be denoted by $\bH^{\infty}_d$ and is the noncommutative analog of the algebra of bounded analytic functions on the unit disc. The algebra $\bH^{\infty}_d$ was called the free semigroup algebra by Davidson and Pitts \cite{DavidsonPitts-inv, DavidsonPitts-alg, DavidsonPitts-NPint}. In fact, $\bH^2_d$ is an NC reproducing kernel Hilbert space and the bounded NC functions on the row ball $\fB_d$ are its left multipliers. Moreover, as was first observed by Popescu \cite{Popescu-funcI}, the algebra of bounded analytic functions on the row ball is unitarily equivalent to the free semigroup algebra (see also \cite[Theorem 3.1]{SSS1}). This explains the notation used. 

In \cite[Theorem 4.6]{DavidsonPitts-alg}, Davidson and Pitts show that every automorphism of $\bH^{\infty}_d$ is automatically WOT-continuous. Moreover, by \cite[Theorem 4.1]{DavidsonPitts-alg}, there is an exact sequence
\[
\xymatrix{1 \ar[r] & \qinn \ar[r] & \aut(\bH^{\infty}_d) \ar[r] & \aut(\B_d) \ar[r] & 1.}
\]
Here $\B_d$ stands for the Euclidean open unit ball in $\C^d$ and $\aut(\B_d)$ is the group of holomorphic automorphisms of $\B_d$. Moreover, there is a section of the quotient map, since every automorphism of $\B_d$ induces a unitarily implemented automorphism of $\bH^{\infty}_d$. In fact, all unitarily implemented automorphisms arise this way. The normal subgroup $\qinn$ is the kernel of the quotient map. Davidson and Pitts call the elements of $\qinn$ quasi-inner automorphisms. In \cite{SSS2}, Salomon, Shalit, and the author prove that automorphisms of $\bH^{\infty}_d$ induce NC analytic automorphisms of $\widetilde{\fB_d}$ that are bi-Lipschitz with respect to a certain analog of a pseudo-hyperbolic metric. Additionally, by \cite[Corollary 6.8]{SSS2}, every quasi-inner automorphism must preserve similarity orbits of irreducible points. Combining this result with \cite[Corollary 6.3]{SSS2} and Corollary \ref{cor:only_irreducible}, we have the following proposition.
\begin{prop} \label{prop:qinn-description}
Let $\Phi \in \qinn$ and let $\varphi \colon \widetilde{\fB_d} \to \widetilde{\fB_d}$ be the induced NC analytic biholomorphism. Then there exist an analytic invertible NC function $G$ on $\fB_d$ and constants $C_1,\, C_2 >0$ such that for every $X \in \fB_d$, 
\begin{itemize}
    \item[(i)] $\varphi(X) = G(X)^{-1} X G(X)$;

    \item[(ii)] $\|\Psi_{\varphi(X)}\| \leq C_1$.

    \item[(iii)] For every $f \in \bH^{\infty}_d$, $\sup_{X \in \fB_d} \|G(X)^{-1} f(X) G(X)\| \leq C_2 \|f\|_{\infty}$.
\end{itemize}
\end{prop}
In particular, item (iii) implies that for every word $\alpha$, we have that $$\sup_{X \in \fB_d} \|G(X)^{-1} X^{\alpha} G(X)\| \leq C_2.$$
\begin{thm} \label{thm:qinn_is_inn}
Every $\Phi \in \qinn$ is inner.
\end{thm}
\begin{proof}
    By \cite[Corollary 4.8]{DavidsonPitts-alg}, for every $1 \leq i \leq d$, $\Phi(L_i) - L_i$ belongs to the WOT-closed two-sided ideal generated by the commutators. In particular, $\Phi(L_i) = L_i + g_i$, where the order of vanishing of $g_i$ at $0$ is at least $2$. Since $\Phi$ is a homomorphism, we have that for every word $\alpha$ with $|\alpha| \geq 1$, $\Phi(L^{\alpha}) = L^{\alpha} + g_{\alpha}$, where the order of vanishing of $g_{\alpha}$ at $0$ is at least $|\alpha| + 1$. Therefore, $(L^{\alpha})^* \Phi(L^{\alpha}) = 1 + (L^{\alpha})^* g_{\alpha} \in \bH^{\infty}_d$. Write the Taylor-Taylor expansion of $G$ at $0$ as $G = \sum_{n=0}^{\infty} \sum_{|\alpha|=n} a_{\alpha} z^{\alpha}$ and similarly write $G^{-1} = \sum_{n=0}^{\infty} \sum_{|\alpha|=n} b_{\alpha} z^{\alpha}$. Note that these series converge uniformly and absolutely on compacta levelwise. Normalize $G(0) = 1$. Then, for every point $X \in \fB_d$ and every word $\alpha$,
    \[
    \sum_{n,m = 0}^{\infty} \sum_{\substack{|\beta| = n\\|\gamma| = m}} b_{\beta} a_{\gamma} X^{\beta \alpha \gamma} = G(X)^{-1} X^{\alpha} G(X) = \Phi(L^{\alpha})(X).
    \]
    Fix a word $\omega$ with $|\omega| < |\alpha|$. We are looking for the coefficient of $z^{\omega}$ in the power series of $(L^{\alpha})^* \Phi(L^{\alpha})$. Looking at the above expression the only candidates are elements of the form $z^{\beta \alpha \gamma}$ such that $z^\beta z^\alpha z^\gamma = z^\alpha z^\omega$. Moreover, we must have that $|\beta| + |\gamma| = |\omega|$. This can only happen for non-empty $\beta$ if $z^\alpha$ is the prefix of $z^\beta z^\alpha$. In order to rule this out, we fix $z^{\alpha} = z_1^m z_2$. Then, for every nonempty $\beta$ of length strictly less than $|\alpha|$, if $z^\beta$ is a prefix of $z^\alpha$, then $z^{\beta} = z_1^k$, for $1 \leq k \leq m$. However, $z^{\beta} z^{\alpha} = z_1^{m+k} z_2$. However, the latter word does not have $z^\alpha$ as a prefix. Therefore, write $f_m = L_2^* L_1^{m *} \Phi(L_1^m L_2)$. We know that $\|f_m\| \leq C$ for all $m \in \N$. Additionally, for $0 \leq k \leq m$, the $k$-th homogeneous component of $f_m$ agrees with the $k$-th homogeneous component of $G$. In particular, $f_m(0) = 1$. Let $f_{m_k}$ be a WOT-convergent subsequence of the $f_m$ with limit $f$. Then, $f \neq 0$, since $f(0) = 1$. Moreover, for every word $\alpha$, we have that $\lim_{k \to \infty} \langle f_{m_k}, z^{\alpha} \rangle = \langle f, z^{\alpha} \rangle$. Since the monomials are orthonormal, this inner product recovers the $\alpha$-th Taylor coefficient. However, for $k$ sufficiently large, the $\alpha$-th Taylor coefficient of $f_{m_k}$ is the same as the one of $G$. We conclude that $f = G$. In particular, $G$ is bounded. Applying the same argument with $\Phi$ replaced by $\Phi^{-1}$, we obtain that $G^{-1}$ is also bounded on $\fB_d$.
\end{proof}

\bibliographystyle{abbrv}
\bibliography{sim}

@preamble{
   "\def\cprime{$'$} "
}

@article {AgMcC-global,
    AUTHOR = {Agler, Jim and McCarthy, John E.},
     TITLE = {Global holomorphic functions in several noncommuting
              variables},
   JOURNAL = {Canad. J. Math.},
  FJOURNAL = {Canadian Journal of Mathematics. Journal Canadien de
              Math\'ematiques},
    VOLUME = {67},
      YEAR = {2015},
    NUMBER = {2},
     PAGES = {241--285},
      ISSN = {0008-414X},
   MRCLASS = {32A99 (32A70)},
  MRNUMBER = {3314834},
MRREVIEWER = {Igor Klep},
       DOI = {10.4153/CJM-2014-024-1},
       
}

@article {AgMcC-implicit,
    AUTHOR = {Agler, Jim and McCarthy, John E.},
     TITLE = {The implicit function theorem and free algebraic sets},
   JOURNAL = {Trans. Amer. Math. Soc.},
  FJOURNAL = {Transactions of the American Mathematical Society},
    VOLUME = {368},
      YEAR = {2016},
    NUMBER = {5},
     PAGES = {3157--3175},
      ISSN = {0002-9947},
   MRCLASS = {47J07 (16S10 32A99)},
  MRNUMBER = {3451873},
MRREVIEWER = {Vicen\c tiu D. R\u adulescu},
       DOI = {10.1090/tran/6546},
       
}

@article {AgMcC-freepick,
    AUTHOR = {Agler, Jim and McCarthy, John E.},
     TITLE = {Pick interpolation for free holomorphic functions},
   JOURNAL = {Amer. J. Math.},
  FJOURNAL = {American Journal of Mathematics},
    VOLUME = {137},
      YEAR = {2015},
    NUMBER = {6},
     PAGES = {1685--1701},
      ISSN = {0002-9327},
   MRCLASS = {30E05 (32A99)},
  MRNUMBER = {3432271},
MRREVIEWER = {Ion Sabba Stefanescu},
       DOI = {10.1353/ajm.2015.0042},
       
}

@article {Ball,
    AUTHOR = {Ball, Joseph A.},
     TITLE = {A lifting theorem for operator models of finite rank on
              multiply-connected domains},
   JOURNAL = {J. Operator Theory},
  FJOURNAL = {Journal of Operator Theory},
    VOLUME = {1},
      YEAR = {1979},
    NUMBER = {1},
     PAGES = {3--25},
      ISSN = {0379-4024},
   MRCLASS = {47A45},
  MRNUMBER = {526287},
MRREVIEWER = {Ronald G. Douglas},
}

@article {BMV-ncrkhs,
    AUTHOR = {Ball, Joseph A. and Marx, Gregory and Vinnikov, Victor},
     TITLE = {Noncommutative reproducing kernel {H}ilbert spaces},
   JOURNAL = {J. Funct. Anal.},
  FJOURNAL = {Journal of Functional Analysis},
    VOLUME = {271},
      YEAR = {2016},
    NUMBER = {7},
     PAGES = {1844--1920},
      ISSN = {0022-1236},
   MRCLASS = {47B32 (47A60)},
  MRNUMBER = {3535321},
MRREVIEWER = {Daniel Seco},
       DOI = {10.1016/j.jfa.2016.06.010},
       
}

@incollection {BMV-int,
    AUTHOR = {Ball, Joseph A. and Marx, Gregory and Vinnikov, Victor},
     TITLE = {Interpolation and transfer-function realization for the
              noncommutative {S}chur-{A}gler class},
 BOOKTITLE = {Operator theory in different settings and related
              applications},
    SERIES = {Oper. Theory Adv. Appl.},
    VOLUME = {262},
     PAGES = {23--116},
 PUBLISHER = {Birkh\"{a}user/Springer, Cham},
      YEAR = {2018},
   MRCLASS = {47B32 (47A60 93B28)},
  MRNUMBER = {3792241},
MRREVIEWER = {Ilie Valu\c{s}escu},
}

@article{Bunce,
	Author = {Bunce, John W.},
	Coden = {JFUAAW},
	Doi = {10.1016/0022-1236(84)90098-3},
	Fjournal = {Journal of Functional Analysis},
	Issn = {0022-1236},
	Journal = {J. Funct. Anal.},
	Mrclass = {47A45 (47A65 47B99)},
	Mrnumber = {744917},
	Mrreviewer = {F. Gilfeather},
	Number = {1},
	Pages = {21--30},
	Title = {Models for {$n$}-tuples of noncommuting operators},
	Volume = {57},
	Year = {1984},
}

@article{DavidsonKennedy_choquet,
  title={Noncommutative {C}hoquet theory},
  author={Davidson, Kenneth R and Kennedy, Matthew},
  journal={arXiv preprint arXiv:1905.08436},
  year={2019}
}

@article {DavidsonPitts-inv,
    AUTHOR = {Davidson, Kenneth R. and Pitts, David R.},
     TITLE = {Invariant subspaces and hyper-reflexivity for free semigroup
              algebras},
   JOURNAL = {Proc. London Math. Soc. (3)},
  FJOURNAL = {Proceedings of the London Mathematical Society. Third Series},
    VOLUME = {78},
      YEAR = {1999},
    NUMBER = {2},
     PAGES = {401--430},
      ISSN = {0024-6115},
     CODEN = {PLMTAL},
   MRCLASS = {47A15 (46L05 47L99)},
  MRNUMBER = {1665248},
MRREVIEWER = {Cecelia Laurie},
       DOI = {10.1112/S002461159900180X},
       
}

@article {DavidsonPitts-alg,
    AUTHOR = {Davidson, Kenneth R. and Pitts, David R.},
     TITLE = {The algebraic stnructure of non-commutative analytic {T}oeplitz
              algebras},
   JOURNAL = {Math. Ann.},
  FJOURNAL = {Mathematische Annalen},
    VOLUME = {311},
      YEAR = {1998},
    NUMBER = {2},
     PAGES = {275--303},
      ISSN = {0025-5831},
   MRCLASS = {47L30 (20E05 46L40 47B35)},
  MRNUMBER = {1625750},
MRREVIEWER = {Ra\'ul E. Curto},
       DOI = {10.1007/s002080050188},
       
}

@article {DavidsonPitts-NPint,
    AUTHOR = {Davidson, Kenneth R. and Pitts, David R.},
     TITLE = {Nevanlinna-{P}ick interpolation for non-commutative analytic
              {T}oeplitz algebras},
   JOURNAL = {Integral Equations Operator Theory},
  FJOURNAL = {Integral Equations and Operator Theory},
    VOLUME = {31},
      YEAR = {1998},
    NUMBER = {3},
     PAGES = {321--337},
      ISSN = {0378-620X},
   MRCLASS = {47A57 (46E22 46L99 47B35 47L80)},
  MRNUMBER = {1627901},
MRREVIEWER = {Scott McCullough},
       DOI = {10.1007/BF01195123},
       
}

@article{EffrosWinkler,
  title={Matrix convexity: operator analogues of the bipolar and {H}ahn--{B}anach theorems},
  author={Effros, Edward G and Winkler, Soren},
  journal={Journal of Functional Analysis},
  volume={144},
  number={1},
  pages={117--152},
  year={1997},
  publisher={Elsevier}
}

@article{Frazho,
	Author = {Frazho, Arthur E.},
	Coden = {JFUAAW},
	Doi = {10.1016/0022-1236(82)90057-X},
	Fjournal = {Journal of Functional Analysis},
	Issn = {0022-1236},
	Journal = {J. Funct. Anal.},
	Mrclass = {47A45 (47A20)},
	Mrnumber = {671311},
	Mrreviewer = {Lawrence Williams},
	Number = {1},
	Pages = {1--11},
	Title = {Models for noncommuting operators},
	Volume = {48},
	Year = {1982},
}

@incollection{HKM12,
  title={Free analysis, convexity and {L}{M}{I} domains},
  author={Helton, J William and Klep, Igor and McCullough, Scott},
  booktitle={Mathematical Methods in Systems, Optimization, and Control: Festschrift in Honor of J. William Helton},
  pages={195--219},
  year={2012},
  publisher={Springer}
}

@article {JuryMartin-affil,
    AUTHOR = {Jury, Michael T. and Martin, Robert T. W.},
     TITLE = {Operators affiliated to the free shift on the free {H}ardy
              space},
   JOURNAL = {J. Funct. Anal.},
  FJOURNAL = {Journal of Functional Analysis},
    VOLUME = {277},
      YEAR = {2019},
    NUMBER = {12},
     PAGES = {108285, 39},
      ISSN = {0022-1236},
   MRCLASS = {46L54 (30H10 46E22 47B32)},
  MRNUMBER = {4019090},
MRREVIEWER = {J. E. Pascoe},
       DOI = {10.1016/j.jfa.2019.108285},
       URL = {https://doi.org/10.1016/j.jfa.2019.108285},
}

@article {JMS-bso,
    AUTHOR = {Jury, Michael T. and Martin, Robert T. W. and Shamovich, Eli},
     TITLE = {Blaschke-singular-outer factorization of free non-commutative
              functions},
   JOURNAL = {Adv. Math.},
  FJOURNAL = {Advances in Mathematics},
    VOLUME = {384},
      YEAR = {2021},
     PAGES = {Paper No. 107720, 42},
      ISSN = {0001-8708},
   MRCLASS = {30H05 (30H10 30J10 46L52)},
  MRNUMBER = {4238916},
       DOI = {10.1016/j.aim.2021.107720},
       URL = {https://doi.org/10.1016/j.aim.2021.107720},
}

@article {JMS-rat,
    AUTHOR = {Jury, Michael T. and Martin, Robert T. W. and Shamovich, Eli},
     TITLE = {Non-commutative rational functions in the full {F}ock space},
   JOURNAL = {Trans. Amer. Math. Soc.},
  FJOURNAL = {Transactions of the American Mathematical Society},
    VOLUME = {374},
      YEAR = {2021},
    NUMBER = {9},
     PAGES = {6727--6749},
      ISSN = {0002-9947},
   MRCLASS = {47A10},
  MRNUMBER = {4302175},
       DOI = {10.1090/tran/8418},
       URL = {https://doi.org/10.1090/tran/8418},
}

@article {JMS-clark,
    AUTHOR = {Jury, Michael T. and Martin, Robert T. W. and Shamovich, Eli},
     TITLE = {Noncommutative rational Clark measures},
   JOURNAL = {Canad. J. Math.},
  FJOURNAL = {Canadian Journal of Mathematics. Journal Canadien de
              Math\'{e}matiques},
      YEAR = {to appear},
       DOI = {10.4153/S0008414X22000384},
       URL = {https://doi.org/10.4153/S0008414X22000384},
}

@book {KVV-book,
    AUTHOR = {Kaliuzhnyi-Verbovetskyi, Dmitry S. and Vinnikov, Victor},
     TITLE = {Foundations of free noncommutative function theory},
    SERIES = {Mathematical Surveys and Monographs},
    VOLUME = {199},
 PUBLISHER = {American Mathematical Society, Providence, RI},
      YEAR = {2014},
     PAGES = {vi+183},
       
}

@article {Popescu-funcII,
    AUTHOR = {Popescu, Gelu},
     TITLE = {Free holomorphic functions on the unit ball of {$B(\mathcal{H})^n$}. {II}},
   JOURNAL = {J. Funct. Anal.},
  FJOURNAL = {Journal of Functional Analysis},
    VOLUME = {258},
      YEAR = {2010},
    NUMBER = {5},
     PAGES = {1513--1578},
      ISSN = {0022-1236},
     CODEN = {JFUAAW},
   MRCLASS = {47A13 (46G20 46L52 47A20)},
  MRNUMBER = {2566311 (2011b:47020)},
MRREVIEWER = {Joseph A. Ball},
       DOI = {10.1016/j.jfa.2009.10.014},
       
}

@article {Popescu-funcI,
    AUTHOR = {Popescu, Gelu},
     TITLE = {Free holomorphic functions on the unit ball of {$B(\mathcal{H})^n$}},
   JOURNAL = {J. Funct. Anal.},
  FJOURNAL = {Journal of Functional Analysis},
    VOLUME = {241},
      YEAR = {2006},
    NUMBER = {1},
     PAGES = {268--333},
      ISSN = {0022-1236},
     CODEN = {JFUAAW},
   MRCLASS = {47A13 (46L54 47A20 47A60)},
  MRNUMBER = {2264252 (2007h:47015)},
MRREVIEWER = {Joseph A. Ball},
       DOI = {10.1016/j.jfa.2006.07.004},
       
}

@article {Popescu-aut,
    AUTHOR = {Popescu, Gelu},
     TITLE = {Free holomorphic automorphisms of the unit ball of {$B(\mathcal{H})^n$}},
   JOURNAL = {J. Reine Angew. Math.},
  FJOURNAL = {Journal f\"ur die Reine und Angewandte Mathematik. [Crelle's
              Journal]},
    VOLUME = {638},
      YEAR = {2010},
     PAGES = {119--168},
      ISSN = {0075-4102},
     CODEN = {JRMAA8},
   MRCLASS = {46L40 (47A45 47B38 47L55 47L80)},
  MRNUMBER = {2595338 (2012c:46161)},
MRREVIEWER = {Joseph A. Ball},
       DOI = {10.1515/CRELLE.2010.005},
       
}

@article{Popescu-dilations,
	Author = {Popescu, Gelu},
	Coden = {TAMTAM},
	Doi = {10.2307/2001359},
	Fjournal = {Transactions of the American Mathematical Society},
	Issn = {0002-9947},
	Journal = {Trans. Amer. Math. Soc.},
	Mrclass = {47A20 (47A45)},
	Mrnumber = {972704},
	Mrreviewer = {L{\'a}szl{\'o} K{\'e}rchy},
	Number = {2},
	Pages = {523--536},
	Title = {Isometric dilations for infinite sequences of noncommuting operators},
	Volume = {316},
	Year = {1989},
}

@article{Popescu-func_calc,
	Author = {Popescu, Gelu},
	Doi = {10.1307/mmj/1029005232},
	Fjournal = {The Michigan Mathematical Journal},
	Issn = {0026-2285},
	Journal = {Michigan Math. J.},
	Mrclass = {47A60 (46L89 47A20 47A45)},
	Mrnumber = {1342494},
	Mrreviewer = {Ra{\'u}l E. Curto},
	Number = {2},
	Pages = {345--356},
	Title = {Functional calculus for noncommuting operators},
	Volume = {42},
	Year = {1995},
}

@article {SSS1,
    AUTHOR = {Salomon, Guy and Shalit, Orr M. and Shamovich, Eli},
     TITLE = {Algebras of bounded noncommutative analytic functions on
              subvarieties of the noncommutative unit ball},
   JOURNAL = {Trans. Amer. Math. Soc.},
  FJOURNAL = {Transactions of the American Mathematical Society},
    VOLUME = {370},
      YEAR = {2018},
    NUMBER = {12},
     PAGES = {8639--8690},
      ISSN = {0002-9947},
   MRCLASS = {46L07 (47L25)},
  MRNUMBER = {3864390},
MRREVIEWER = {\'{E}ric Ricard},
       DOI = {10.1090/tran/7308},
       URL = {https://doi.org/10.1090/tran/7308},
}

@article {SSS2,
    AUTHOR = {Salomon, Guy and Shalit, Orr M. and Shamovich, Eli},
     TITLE = {Algebras of noncommutative functions on subvarieties of the
              noncommutative ball: the bounded and completely bounded
              isomorphism problem},
   JOURNAL = {J. Funct. Anal.},
  FJOURNAL = {Journal of Functional Analysis},
    VOLUME = {278},
      YEAR = {2020},
    NUMBER = {7},
     PAGES = {108427, 54},
      ISSN = {0022-1236},
   MRCLASS = {46L52 (46L07 47L25 47L80)},
  MRNUMBER = {4053626},
MRREVIEWER = {Volker Runde},
       DOI = {10.1016/j.jfa.2019.108427},
       URL = {https://doi.org/10.1016/j.jfa.2019.108427},
}

@book {AMcCY-book,
    AUTHOR = {Agler, Jim and McCarthy, John Edward and Young, Nicholas},
     TITLE = {Operator analysis---{H}ilbert space methods in complex
              analysis},
    SERIES = {Cambridge Tracts in Mathematics},
    VOLUME = {219},
 PUBLISHER = {Cambridge University Press, Cambridge},
      YEAR = {2020},
     PAGES = {xv+375},
}

@article {GMS-boundaries,
    AUTHOR = {Griesenauer, Erin and Muhly, Paul S. and Solel, Baruch},
     TITLE = {Boundaries, bundles, and trace algebras},
   JOURNAL = {New York J. Math.},
  FJOURNAL = {New York Journal of Mathematics},
    VOLUME = {24A},
      YEAR = {2018},
}

@article {Artin-azumaya,
    AUTHOR = {Artin, M.},
     TITLE = {On {A}zumaya algebras and finite dimensional representations
              of rings},
   JOURNAL = {J. Algebra},
  FJOURNAL = {Journal of Algebra},
    VOLUME = {11},
      YEAR = {1969},
     PAGES = {532--563},
}

@article {Procesi-invariants,
    AUTHOR = {Procesi, C.},
     TITLE = {The invariant theory of {$n\times n$} matrices},
   JOURNAL = {Advances in Math.},
  FJOURNAL = {Advances in Mathematics},
    VOLUME = {19},
      YEAR = {1976},
    NUMBER = {3},
     PAGES = {306--381},
}

@article {SerKha-domains,
    AUTHOR = {Sergeev, A. G. and Kha{\u\i}ntsner, P.},
     TITLE = {The extended matrix disk is a domain of holomorphy},
   JOURNAL = {Izv. Akad. Nauk SSSR Ser. Mat.},
  FJOURNAL = {Izvestiya Akademii Nauk SSSR. Seriya Matematicheskaya},
    VOLUME = {55},
      YEAR = {1991},
    NUMBER = {3},
     PAGES = {647--657},
}

@book {Bhatia-book,
    AUTHOR = {Bhatia, Rajendra},
     TITLE = {Matrix analysis},
    SERIES = {Graduate Texts in Mathematics},
    VOLUME = {169},
 PUBLISHER = {Springer-Verlag, New York},
      YEAR = {1997},
     PAGES = {xii+347},
}

@book {Taylor-sev_comp_var,
    AUTHOR = {Taylor, Joseph L.},
     TITLE = {Several complex variables with connections to algebraic
              geometry and {L}ie groups},
    SERIES = {Graduate Studies in Mathematics},
    VOLUME = {46},
 PUBLISHER = {American Mathematical Society, Providence, RI},
      YEAR = {2002},
     PAGES = {xvi+507},
}

@article {Frenkel-nonab_coh,
    AUTHOR = {Frenkel, Jean},
     TITLE = {Cohomologie non ab\'elienne et espaces fibr\'es},
   JOURNAL = {Bull. Soc. Math. France},
  FJOURNAL = {Bulletin de la Soci\'et\'e{} Math\'ematique de France},
    VOLUME = {85},
      YEAR = {1957},
     PAGES = {135--220},
}

@book {Hatcher-book,
    AUTHOR = {Hatcher, Allen},
     TITLE = {Algebraic topology},
 PUBLISHER = {Cambridge University Press, Cambridge},
      YEAR = {2002},
     PAGES = {xii+544},
}

@book {HilNeeb-lie_groups,
    AUTHOR = {Hilgert, Joachim and Neeb, Karl-Hermann},
     TITLE = {Structure and geometry of {L}ie groups},
    SERIES = {Springer Monographs in Mathematics},
 PUBLISHER = {Springer, New York},
      YEAR = {2012},
     PAGES = {x+744},
}

@book {BanSta-alg_methods,
    AUTHOR = {B\u{a}nic\u{a}, Constantin and St\u{a}n\u{a}\c{s}il\u{a}, Octavian},
     TITLE = {Algebraic methods in the global theory of complex spaces},
      NOTE = {Translated from the Romanian},
 PUBLISHER = {Editura Academiei, Bucharest; John Wiley \& Sons, London-New
              York-Sydney},
      YEAR = {1976},
     PAGES = {296},
}

@book {SGA1,
    AUTHOR = {Grothendieck, Alexander and Raynaud, M.},
     TITLE = {Rev\^etements \'etales et groupe fondamental},
      NOTE = {S\'eminaire de G\'eom\'etrie Alg\'ebrique du Bois Marie, 1960/61},
 PUBLISHER = {ISpringer Berlin, Heidelberg},
      YEAR = {1971},
}

@article {HochRob-CM,
    AUTHOR = {Hochster, Melvin and Roberts, Joel L.},
     TITLE = {Rings of invariants of reductive groups acting on regular
              rings are {C}ohen-{M}acaulay},
   JOURNAL = {Advances in Math.},
  FJOURNAL = {Advances in Mathematics},
    VOLUME = {13},
      YEAR = {1974},
     PAGES = {115--175},
}

@book {Hirsch-diff_top,
    AUTHOR = {Hirsch, Morris W.},
     TITLE = {Differential topology},
    SERIES = {Graduate Texts in Mathematics},
    VOLUME = {33},
      NOTE = {Corrected reprint of the 1976 original},
 PUBLISHER = {Springer-Verlag, New York},
      YEAR = {1994},
     PAGES = {x+222},
}

@incollection {LeBruynProcesi-etale,
    AUTHOR = {Le Bruyn, Lieven and Procesi, Claudio},
     TITLE = {\'Etale local structure of matrix invariants and concomitants},
 BOOKTITLE = {Algebraic groups {U}trecht 1986},
    SERIES = {Lecture Notes in Math.},
    VOLUME = {1271},
     PAGES = {143--175},
 PUBLISHER = {Springer, Berlin},
      YEAR = {1987},
}
\vspace{0.5cm}
\end{document}